\documentclass[submission]{dmtcs}

\usepackage[latin1]{inputenc}
\usepackage{subfigure}
\usepackage[square]{natbib}

\usepackage{
amsfonts,amsmath,amssymb}

\newcommand{\eq}{\begin{equation}}
\newcommand{\en}{\end{equation}}

\newcommand{\giv}{\,|\,}

\newcommand{\prob}{\mathbb P}
\newcommand{\ex}{\mathbb E}

\newcommand{\Nat}{\Bbb N}

\newtheorem{theorem}{\large Theorem}%[section]
\newtheorem{proposition}[theorem] {\large Proposition}
\newtheorem{definition}[theorem]{\large Definition}

\RCSdef$Revision: 1.3 $\endRCSdef
\rcsMajMin

\author{Alexander Gnedin\addressmark{1}\thanks{gnedin@math.uu.nl}
}
\title[Constrained  exchangeable partitions]{ Constrained  exchangeable partitions}
\address{\addressmark{1} Mathematical Institute, Utrecht University, P.O. Box 80010, 3508 TA Utrecht, The Netherlands 
}

\received{2 April 2006}
\revised{20 June 2006}
\accepted{27 June 2006}

\begin{document}

\maketitle

\begin{abstract}
 % \begin{quotation}
 %   This is revision {\rcsMaj.\rcsMin} of this document.
 % \end{quotation} 
 \begin{description}
For a class of random partitions of an infinite 
set a de Finetti-type representation is derived,
and in one special case a central limit theorem for the number of blocks is shown.
  \end{description}
\end{abstract}

%\clearpage
%\tableofcontents

\section{Introduction}\label{intro}

Under a {\it partition} of the set $\Nat$ we shall mean
a sequence  $(b_1,b_2,\ldots)$ of  subsets of $\Nat$ such that (i) the sets $b_j$ are disjoint, (ii) $\cup_j b_j=\Nat$,
(iii) if $b_k=\varnothing$ then also $b_{k+1}=\varnothing$ and (iv) if $b_{k+1}\neq\varnothing$ then 
$\min b_k<\min b_{k+1}$.  
Condition
(iv) says that the sequence of minimal elements of the blocks is increasing.
One can  think of partition as a mapping which sends a generic  element $j\in \Nat$ to one of the infinitely many blocks,
in such a way that conditions (iii) and (iv) are fulfilled.

\par A random partition $\Pi=(B_k)$ of $\Nat$ (so, with random blocks $B_k$) is a random variable with values in the set of partitions of $\Nat$.
This concept can be made precise by means of a projective limit construction and 
 the measure extension theorem. To this end, one identifies $\Pi$ with  consistent
partitions $\Pi_n:=\Pi|_{[n]}$ ($n=1,2,\ldots$) of finite sets $[n]:=\{1,\ldots,n\}$.
Note that the restriction $\Pi_n$, which is obtained by removing all elements not in $[n]$,
still has the blocks in the order of increase of their least elements.

\par There is a well developed theory of exchangeable partitions \cite{BertoinBook, Kallenberg,  CSP}.
 Recall that  $\Pi=(B_j)$ 
is {\it exchangeable} if the law of $\Pi$ is invariant under all
bijections $\sigma:\Nat\to\Nat$. 
Partitions with weaker symmetry properties have also been studied.
Pitman \cite{PTRF} introduced {\it partially exchangeable} random partitions of $\Nat$ with  the property that 
the law of $\Pi$ is invariant under all bijections $\sigma:\Nat\to\Nat$ that preserve the order of blocks,
meaning that the sequence of
the least elements of the sets $\sigma(B_1),\sigma(B_2),\ldots$ is also increasing.
Pitman \cite{PTRF}
derived a de Finetti-type representation for partially exchangeable partitions
and established a criterion for their exchangeability.
Kerov \cite{KerovSub} studied a closely related  structure of virtual permutations of $\Nat$,
which may be seen as partially exchangeable partitions with some total ordering of elements within
each of the blocks.
Kallenberg \cite{Kallenberg} characterised {\it spreadable} partitions whose law is invariant under increasing 
injections $\sigma:\Nat\to\Nat$.

\par In this note we consider constrained random partitions of $\Nat$ which satisfy the condition that, 
for a fixed integer sequence $\rho=(\rho_1,\rho_2,\ldots)$ with $\rho_k\geq 1$,
each block   $B_k$ contains
$\rho_k$ least elements 
of $\cup_{j\geq k}B_j$, for every $k$ with $B_k\neq\varnothing$.
It is easy to check that this condition holds if and anly if the sequence comprised 
of $\rho_1$ least elements of
$B_1$, followed by  $\rho_2$ least elements of $B_2$ and so on, is itself 
an increasing sequence. We shall focus on the constrained partitions with 
the following 
symmetry property.
\begin{definition}
{\rm For a given sequence $\rho$,
we call $\Pi$ {\it constrained exchangeable} if $\Pi$ 
is a constrained partition with respect to $\rho$ and the law 
of $\Pi$
is invariant under all
bijections $\sigma:\Nat\to\Nat$ that preserve this property. 
}
\end{definition}

\noindent
Since the law of $\Pi$ is uniquely determined by the laws of finite restrictions $\Pi_n$,
the constrained exchangeability of $\Pi$ amounts to the analogous property of $\Pi_n$'s for each $n=1,2,\ldots$.
To gain a feeling of the property, the reader is suggested to check that for $\rho=(1,2,1,\ldots)$
the partition $\Pi_8$ assumes
the values $(\{1,3,5\},\{2,4,6\},\{7,8\})$ and
$(\{1,2,3\},\{4,5,8\},\{6,7\})$ with the same probability.

\par Every  partition of $\Nat$ is constrained
with respect to $\rho=(1,1,\ldots)$,
and every constrained exchangeable partition with this $\rho$ is
partially exchangeable in the sense of  Pitman \cite{PTRF}.
In principle, any constrained exchangeable partition may be reduced to some Pitman's partially exchangeable 
partition by isolating $\rho_k-1$ least elements of $B_k$ in $\rho_k-1$ singleton blocks, for each $\rho_k>1$, but
this viewpoint will not be adopted here.
\par  For general  
$\rho\neq(1,1,\ldots)$ the constrained exchangeable partitions which
are also exchangeable are rather uninteresting, since they 
cannot have infinitely many blocks:

\begin{proposition}
Let $\Pi$ be a constrained partition 
with  respect to some $\rho$ which has  $\rho_k>1$ for some  $k$.
If~
$\Pi$ is exchangeable then $\Pi$ has 
at most $k$ nonempty blocks.
\end{proposition}
\begin{proof}
Suppose $B_k\neq\varnothing$, then by Kingman's representation of exchangeable partitions \cite{CSP} the set
$\cup_{j\geq k}B_j$ contains infinitely many elements.
For the same reason  $\#B_k\geq 2$ implies that $B_k$ is an infinite
set, and that partition $\Pi'$ obtained by restricting $\Pi$ to $\cup_{j\geq k}B_j$ and re-labelling  
the elements of $\cup_{j\geq k}B_j$  by $\Nat$ in increasing order is an exchangeable partition of $\Nat$.
But then with probability one $\Pi'$ is the trivial single-block partition, 
because elements $1$ and $2$ are always 
in the same block. 
\end{proof}

\par In many contexts where random partitions appear, exchangeability is an obvious kind of symmetry.
Constrained exchangeability may appear
 when some initial elements of the blocks play 
a special role of `establishing' the block.
 To illustrate, consider the following situation.
Suppose there is a sequence of independent random points sampled from  some distribution on ${\mathbb R}^d$.
Define $D_1$ as the convex hull of the first $\rho_1$ points,
$D_2$ as the convex hull of the first $\rho_2$ points not in $D_1$, 
$D_3$ as the convex hull of the first $\rho_3$ points not in $D_1\cup D_2$, 
etc.
Divide ${\mathbb R}^d$ in disjoint nonempty subsets $G_1=D_1, G_2= D_2\setminus D_1, G_3=D_3\setminus(D_1\cup D_2),\ldots$.
A constrained exchangeable partition $\Pi$ 
of $\Nat$ 
is defined then by assigning to block $B_k$ the indices of $\rho_k$ 
initial points that determine  $D_k$ and the indices of all further points  that 
hit $G_k$. 
\par Of course,
there is nothing special in 
the convex hulls construction, and any other  way of  `spanning' a spatial domain $D_k$ on $\rho_k$ sample
points and then `peeling' the space in $G_k$'s will also result in a constrained exchangeable partition.
An example of this kind related to 
 multidimensional
records will be given.

\par In what follows we extend Pitman's \cite{PTRF} sequential realisation  of partitions
 via frequencies of blocks,
to cover arbitrary constrained exchangeable partitions.
Generalising a result on exchangeable partitions
\cite{Sieve}
we shall also derive a central limit theorem  for the number of blocks of finite partitions 
$\Pi_n=\Pi|_{[n]}$ in one important case of  partitions induced by  a `stick-breaking' scheme.

\section{Constrained sampling}
We fix throughout 
a sequence of positive integers $\rho$.
Recall that a  {\it composition} is a finite sequence of positive integers called parts,
e.g. $(3,1,2)$ is a composition of $6=3+1+2$ with three parts.
We say that a composition
$\lambda=(\lambda_1,\ldots,\lambda_\ell)$ 
is a {\it constrained composition of} $n$ if 
$\lambda_j\geq\rho_j$ for $j=1,\ldots,\ell-1$ and $|\lambda|:=\sum\lambda_j=n$.
\par For each $\lambda$ a constrained composition of $n$, the following random algorithm, which may be called
{\it constrained sampling}, yields another constrained composition $\mu$ of $n-1$.
Imagine a row of  boxes labeled $1,\ldots,\ell$ and occupied by  $\lambda_1,\ldots,\lambda_\ell$  white balls. 
Let $\Lambda_j:=\lambda_j+\ldots+\lambda_\ell\,,\,j\leq \ell$.
At the first step,
$\rho_1$ balls in box 1 are re-painted black and then
 a white ball is drawn uniformly at random from all $\Lambda_1-\rho_1$ white balls. 
If the ball drawn was in box 1, the ball is deleted and the
new composition is $\mu=(\lambda_1-1,\lambda_2,\ldots,\lambda_\ell)$, and
if the ball drawn was in some other box, it is returned to the box and the process continues,
so that 
at the second step
$\rho_2$ balls in box 2 are re-painted black and a white ball is drawn uniformly at random from
boxes $2,\ldots,\ell$. If the  second ball drawn was in box 2, the ball is  deleted and  the new composition is
$\mu=(\lambda_1,\lambda_2-1,\lambda_3,\ldots,\lambda_\ell)$, and so on.
If the procedure does not terminate in $\ell-1$ steps, then
a ball is deleted from the last box and the new composition is
 $\mu=(\lambda_1,\ldots,\lambda_{\ell-1},\lambda_\ell-1)$.
By this description, for $j<\ell$ the transition probability from $\lambda$ to $\mu=(\ldots,\lambda_j-1,\ldots)$
is 
$${\Lambda_2\over (\Lambda_1-\rho_1)}\cdots{\Lambda_{j}\over (\Lambda_{j-1}-\rho_{j-1})}
{(\lambda_j-\rho_j)\over(\Lambda_j-\rho_j)}\,,$$
while the transition probability from $\lambda$ to $\mu=(\cdots,\lambda_\ell-1)$ is
$${\Lambda_2\over (\Lambda_1-\rho_1)}\cdots{\Lambda_{\ell}\over (\Lambda_{\ell-1}-\rho_{\ell-1})}\,.$$
\par A random sequence ${\cal C}=({\cal C}_n)$ of constrained compositions of integers $n=1,2,\ldots$ is 
called  {\it consistent}
if  ${\cal C}_{n-1}$
has the same law as the composition derived from ${\cal C}_n$ by 
the above constrained sampling procedure, for each $n>1$. 
Every consistent sequence $({\cal C}_n)$ is an inverse Markov chain
with some co-transition probabilities  depending only on $\rho$.
By analogy with \cite{gnedin97} a consistent sequence $\cal C$ will be called a 
{\it constrained composition structure}. 

\par 
If the constraints are determined by  $\rho=(1,1,\ldots)$, the constrained sampling amounts
to a co-transition rule related to the  partially exchangeable partitions in \cite{PTRF}. 
The unconstrained sampling (corresponding to $\rho=(0,0,\ldots)$) leads to composition structures  studied in \cite{gnedin97, RCS, RPS}.

\section{Basic representation}\label{basic}
For $\Pi$ a constrained partition of $\Nat$ with blocks $B_1,B_2,\ldots$
we define,
for each $n=1,2,\ldots$, a composition
${\cal C}_n$ of $n$ as the  finite sequence of positive values in  
$\#(B_1\cap[n]),\#(B_2\cap[n]),\ldots$. 
We call this composition the {\it shape} of $\Pi_n$ and write ${\cal C}_n={\tt shape}(\Pi_n)$.

\par The number of constrained  partitions of $[n]$ with shape $\lambda$ is equal to
\eq\label{d}
d(\lambda):=\prod_{j=1}^{\ell -1} {\Lambda_j-\rho_j\choose \lambda_j-\rho_j}\,.
\en
Similarly, the number of partitions of $[n]$ with    shape
$\lambda=(\lambda_1,\ldots,\lambda_\ell)$ and whose restriction on $[n']$ (for $n'<n$) has shape
$\mu=(\mu_1,\ldots,\mu_k)$  
is equal to
\eq\label{dd}
d(\lambda,\mu):=\left[\prod_{j=1}^{\ell-1}{{\rm M}_j-\Lambda_j\choose \mu_j-\lambda_j}\right]
{{\rm M}_\ell-\Lambda_\ell-(\rho_\ell-\lambda_\ell)_+\choose \mu_\ell-\rho_\ell\vee\lambda_\ell}
\left[\prod_{j=\ell+1}^{k-1} {{\rm M}_j-\rho_j\choose \mu_j-\rho_j}\right],
\en
where ${\rm M}_j=\mu_j+\ldots+\mu_k\,,\,j\leq k$.

\par Introduce a function of compositions
$$ p(\lambda):=\prob({\tt shape}(\Pi_n)=\lambda).$$
It is easy to check that the consistency of $\Pi_n$'s with respect to restrictions
implies that the ${\cal C}_n$'s are consistent in the sense of constrained sampling,
therefore appealing to Kolmogorov's measure extension theorem we have:

\begin{proposition} 
The formula  
$$\prob(\Pi_n=\,\cdot\,) = p({\tt shape}(\,\cdot\,))/d({\tt shape}(\,\cdot\,))$$ 
establishes a 
 canonical homeomorphism 
between the distributions of
constrained exchangeable partitions of $\Nat$
and constrained composition structures. Conditionally given ${\cal C}_n={\tt shape}( \Pi_n)=\lambda$
the distribution of $\Pi_n$ is uniform on the set of constrained partitions of $[n]$ with shape $\lambda$. 
\end{proposition}

\par The following basic construction modifies  the one exploited   in \cite{KerovSub, PTRF}. 
Let $(P_1,P_2,\ldots)$ be an arbitrary sequence of random variables
satisfying $P_k\geq 0$ and $\sum_k P_k\leq 1$. 
A constrained exchangeable partition $\Pi$ directed by $(P_k)$
is defined as follows. 
Conditionally given $(P_k)$ the partition 
is obtained by successive extension of $\Pi_n$ to $\Pi_{n+1}$, 
for each $n=1,2,\ldots$, according to the rules:
given $\Pi_n$ with ${\tt shape}(\Pi_n)=(\lambda_1,\ldots,\lambda_\ell)$, the element $n+1$
\begin{itemize}
\item[{\rm (i)}] joins the block $B_j,~ j<\ell,$ with probability $P_j$,
\item[{\rm (ii)}] if $\lambda_\ell<\rho_\ell$ joins the  block $B_\ell$ with probability
$1-\sum_{j=1}^{\ell-1}P_j$,
\item[{\rm (iii)}] and if $\lambda_\ell\geq \rho_\ell$ joins the  block $B_\ell$
with probability $P_\ell$ or starts the new block  $B_{\ell+1}$ with probability $1-\sum_{j=1}^\ell P_j$.
\end{itemize}
Explicitly, for the  function
$p$ of compositions
we have the formula
\eq\label{Fin}
p(\lambda)
=d(\lambda)\,{\mathbb E}\,\left\{
\left[ \prod_{j=1}^{\ell-1} \left(1-\sum_{i=1}^{j-1}P_i\right)^{\rho_j}P_j^{\lambda_j-\rho_j}
\right]
\left(1-\sum_{i=1}^{\ell-1}P_i
\right)^{\rho_\ell\wedge\lambda_\ell} P_\ell^{(\lambda_\ell-\rho_\ell)_+}\right\}\,.
\en

\par The next de Finetti-type result states that the construction covers all possible constrained exchangeable partitions.
The proof is only sketched, since it follows the same lines as in 
\cite{KerovSub, PTRF}.

\begin{proposition} For $\Pi$ a constrained exchangeable partition, the 
normalised shapes ${\tt shape}(\Pi_n)/n$ (considered as sequences padded by infinitely many zeroes)
converge in the product topology with probability one
to some random limit $(P_1,P_2,\ldots)$ satisfying $P_j\geq 0$ and  $\,\sum_j P_j\leq 1$.
Conditionally given $(P_k)$ the partition $\Pi$ is recovered according to the above rules {\rm (i)-(iii)}.
\end{proposition}
\begin{proof} The key point is to show the existence of frequencies. This can be concluded
from de Finetti's theorem for $0-1$ exchangeable sequences by noting that the indicators 
$1(m~{\rm belongs~to~block}~B_k)$ 
for $m>n$ are conditionally exchangeable
given that the block $B_k$ has at least $\rho_k$ representatives in $[n]$. Alternatively,
one can use, as in \cite{KerovSub}, more direct Martin boundary arguments which exploit the
explicit formulas (\ref{d}) and (\ref{dd}) to show that 
the pointwise limit of ratios
${d(\,\cdot\,,\mu)/d(\mu)}$, as $m=|\mu|\to\infty$, exists if and only if $\mu_j/m$ converge for every
$j$.
\end{proof}

\section{The formation sequence}

Nacu \cite{Nacu} established that the law of a partially exchangeable partition is uniquely determined by the law
of the increasing sequence of the least elements of blocks.
We show that a similar result holds for every constrained exchangeable partition $\Pi$ with 
general $\rho$.

\par We define the {\it formation sequence} to be the sequence obtained by selecting the $\rho_k$th  least 
element of the block $B_k$ of $\Pi$, for  $k=1,2,\ldots$.
For a composition $\lambda$ let $q(\lambda)$ be the probability that the formation sequence
starts with elements
$\lambda_1,\lambda_1+\lambda_2,\ldots,\lambda_1+\ldots+\lambda_\ell$.
Let  $(P_j)$ be the frequencies as in Section \ref{basic}, and
introduce the variables 
$$
H_k=1-\sum_{j=1}^{k}P_j\,,~~{\rm so~~}P_k=H_{k-1}-H_k\,,
$$ 
where we set $H_0=1$.
Then
\eq\label{qf}
q(\lambda_1,\ldots,\lambda_\ell)=
{\mathbb E}
\left[\prod_{j=1}^{\ell-1}
{\lambda_{j+1}-1\choose \rho_{j+1}-1}H_j^{\rho_{j+1}}
(1-H_j)^{\lambda_{j+1}-\rho_{j+1}}\right].
\en
Comparing this with (\ref{Fin}) written in the same variables (where $H_0=1$) we obtain
for constrained compositions
\eq\label{pf}
p(\lambda)=d(\lambda)\,{\mathbb E}\,\left\{
\left[ \prod_{j=1}^{\ell-1} H_{j-1}^{\rho_j}(H_{j-1}-H_j)^{\lambda_j-\rho_j}
\right]
H_{\ell-1}^{\rho_\ell\wedge\lambda_\ell}(H_{\ell-1}-H_\ell)^{{(\lambda_\ell-\rho_\ell)_+}}\right\},
\en
\noindent
which leads to the following conclusion:
\begin{proposition} There is an invertible linear transition from 
$p$ to $q$. Hence each of these two functions on compositions uniquely determines the law of $\Pi$.
\end{proposition}
  \begin{proof} 
    The substantial part of the claim is showing that we can compute $p$ from $q$.
     To that end, start by observing that $p$ is uniquely determined by the values on 
     compositions of the type $(\lambda_1,\ldots,\lambda_{\ell-1},\rho_\ell)$. To see that this
     follows from the consistency for various $n$, argue by
    induction in $m=0,\ldots,\rho_\ell$ for compositions $(\,\ldots,\,\rho_\ell-m)$. 
    Now, for such compositions whose last part meets the constraint exactly,
    (\ref{pf}) and (\ref{qf}) involve the same factors of the type $H_j^{\rho_j}$,
    hence $p$ can be reduced to $q$ by expanding each factor
    $(H_{j-1}-H_j)^k=((1-H_j)-(1-H_{j-1}))^k$ using the binomial formula.
    \end{proof}

\section{The paintbox}

Paintbox representations based on the uniform sampling from $[0,1]$
are often  used to model  exchangeable structures and their relatives
\cite{ gnedin97, GO, RCS, CSP}. We shall design  a version that  is appropriate  for  constrained exchangeable 
partitions.

\par Let $1=H_0\geq H_1\geq H_2\geq\ldots\geq 0$ be an arbitrary nonincreasing random sequence.
Let $(U_n)$ be a sequence of independent $[0,1]$-uniform random points, independent of 
$(H_k)$. We define a new sequence $(U_n)\giv_{\rho}(H_k)$ with some of $U_n$'s replaced by 
$H_k$'s, as follows.
Replace $U_1,\ldots ,U_{\rho_1},$ by $H_1$. 
Then replace the first $\rho_2$ entries which belong to $U_{\rho_1},U_{\rho_1+1},\ldots$ and hit 
$[0, H_1[$ by $H_2$.
Inductively, when  $H_1,\ldots,H_k$ get used, respectively, $\rho_1,\ldots,\rho_k$ times, 
keep on screening uniforms until replacing the first $\rho_{k+1}$ points hitting $[0,H_k[\,$
by $H_{k+1}$. 
Eventually all $H_k$'s will enter the resulting sequence.
\par The construction has an interpretation in terms of the classical theory of records
(see \cite{chain, KerovSub} for a special case).

\begin{proposition}
Conditionally given $(H_k)$,
the sequence $(U_n)\giv_{\rho}(H_k)$ has the same distribution as $(U_n)$ 
conditioned on 
the event that the sequence of lower records in $(U_n)$ is $(H_k)$, with 
the record value $H_k$ repeated $\rho_k$ times. 
\end{proposition}

In this framework, we define a partition $\Pi$
by assigning to block $B_k$ all integers which label 
the entries   of $(U_n)\giv_{\rho}(H_k)$ falling in 
$[H_{k},H_{k-1}[\,$. 
Given $(H_k)$, the chance for
$U_n$ to hit $[H_j,H_{j-1}[\,$ is $P_j=H_{j-1}-H_j$, therefore the construction is equivalent to that defined by the rules
(i)-(iii) above.

\section{Stick-breaking partitions}\label{6}

Explicit evaluation of the function $p$ is possible when the frequencies involve a kind of
independence. To this end,
it is convenient to 
introduce yet another set of  variables $(W_k)$ (sometimes called {\it residual fractions}) which satisfy $H_k=W_1\cdots W_k$,
$W_k\in [0,1]$. In these variables (\ref{pf}) becomes

\eq\label{wf}
p(\lambda)=d(\lambda)\,{\mathbb E}\left\{\left[\prod_{j=1}^{\ell-1}W_j^{\Lambda_j-\lambda_j}(1-W_j)^{\lambda_j-\rho_j}\right]
(1-W_\ell)^{(\lambda_\ell-\rho_\ell)_+}\right\}
\en
where $\Lambda_j=\lambda_j+\ldots+\lambda_\ell$. As in \cite{PTRF}, if the $W_k$'s 
are independent, (\ref{wf}) assumes the product form  
\eq\label{prod}
p(\lambda)=\prod_{k=1}^\ell q_k(\Lambda_k:\lambda_k)
\en
with the {\it decrement matrices}
\eq\label{dec}
q_k(n:m)={n-\rho_k\choose m-\rho_k} \ex \left[ (1-W_k)^{m-\rho_k}W_k^{n-m}\right]\,,~~~1\leq m\leq n,
\en
and the convention ${-i\choose -j}=1(i=j)$ for negative arguments of the binomial coefficients. 
In fact, 
(\ref{prod}) forces representation (\ref{dec})
(this fact is implicit in \cite{KerovSub, PTRF} in the case of partially exchangeable partitions):

\begin{proposition} A constrained composition structure $({\cal C}_n)$
satisfies {\rm (\ref{prod})}
with some decrement matrices $q_k\,,\,k=1,2,\ldots$, if and only if there exist independent $[0,1]$-valued random variables
$(W_k)$ such that {\rm (\ref{dec})} holds.
\end{proposition}
\begin{proof} We argue the `only if' part.
For $n<\rho_1$ we have $q_1(n:n)=1$ by definition. For $n\geq \rho_1$ 
the constrained sampling consistency yields
$$q_1(n:m)={m+1-\rho_1\over n+1-\rho_1}\,q_1(n+1:m+1)+
{n+1-m\over n+1-\rho_1}\,q_1(n+1:m)
$$
which is the familiar Pascal-triangle recursion in the variables $n-\rho_1,~m-\rho_1$,
therefore the integral representation (\ref{dec}) follows as a known consequence of the Hausdorff moments problem.
The case $k>1$ is completely analogous. The independence of the $W_k$'s is obvious from (\ref{prod}).
\end{proof}

\noindent
We note in passing that the product formula (\ref{prod}) with a single decrement
matrix leads, in a related setting of regenerative composition structures,  to a nonlinear recursion and a very different conclusion
\cite{RCS}. See \cite{Gibbs} for product formulas of  another kind 
in the exchangeable case.

\par Suppose now that  $W_k$'s are independent and have  beta$(a_k,b_k)$ distributions, 
whose density is
$$(1-s)^{a_k-1}s^{b_k-1}/{\rm B}(a_k,b_k).$$ 
The  rows of the decrement matrices 
are then P{\'o}lya-Eggenberger distributions
$$
q_k(n:m)={n-\rho_k\choose m-\rho_k}{(a_k)_{m-\rho_k}(b_k)_{n-m}\over (a_k+b_k)_{n-\rho_k}}\,,~~~~~m=1,\ldots,n.
$$
For instance, taking positive integer $a_k,\,b_k$ and $\rho_k=a_k+b_k-1$, a partition
$\Pi$ is constructed as follows:
replace $U_1,\ldots,U_{\rho_1}$ by the value $H_1$ equal to
 the $b_1$th minimal  order statistic of these points,
then replace the first $\rho_2$ uniforms that hit $[0,H_1[$ by the value $H_2$ equal to the $b_2$th minimal order statistic of 
these hits, etc, thus defining a partition via $(U_n)\giv_\rho(H_k)$. 
 A distinguished class of structures of this kind is the Ewens-Pitman two-parameter family of exchangeable 
partitions \cite{ selfsim, PTRF, CSP} with $\rho=(1,1,\ldots)$,
$a_k=\theta+k\alpha$ and $b_k=1-\alpha$ (for suitable $\alpha$ and $\theta$); 
the product formula simpifies in this case due to a major telescoping of factors.

\section{Counting the blocks}\label{s:CLT}
Let $K_n$ be the number of blocks of $\Pi_n$,
 which in the $(U_j)\giv_\rho (H_k)$-representation coincides with
the number of intervals $]H_1,H_0], \,]H_2,H_1],\ldots$ discovered by the $n$ first terms of the sequence.
Conditionally given $(H_k)$, $K_n$ is the  
number of certain
independent geometric summands which sum to no more than $n$. In particular,
the difference between the $k$th and the $(k+1)$st  entries of the formation sequence
follows the negative binomial distribution with parameters $\rho_k, H_k$.
\par We shall proceed by assuming a `stick-breaking' scheme $H_k=W_1\cdots W_k\,,~k=1,2,\ldots,$ with independent identically 
distributed $W_k$'s.
We assume further that the logarithmic 
moments $\mu={\mathbb E}[-\log W_1],$ $\sigma^2={\rm Var}[-\log W_1]$ are both finite.
The idea is to derive a CLT for $K_n$ from the standard CLT for renewal processes \cite{Feller}.
Similar technique was used in \cite{Sieve, chain}, but in 
the new stuation  we need to also limit the growth of $\rho_k$  as $k\to\infty$.

\par We will show that $K_n$ is asymptotic to $J_n:=\max\{k: H_k>1/n\}$. The last quantity
is indeed asymptotically Gaussian with the mean $(\log n)/\mu$ and the variance $(\log n)/(\sigma^2\mu^{-3})$
because $J_n$ is just the 
number of renewal epochs within  $[0,\log n]$ of the renewal process with steps $-\log W_k$. 
In loose terms, we will exploit a `cut-off phenomenon': typically, only a few  points 
out of $n$ uniforms fall  below
$1/n$, while for $k<J_n$ essentially all intervals get hit, with exponentially growing occupancy numbers
when scanned backwards in  $k$ from  $k=J_n$.

\begin{proposition}\label{CLT} 
Suppose $\Pi$ is directed by $H_k=W_1\cdots W_k,$ where for $k=1,2,\ldots$ the $W_k$'s 
are i.i.d. with  finite logarithmic moments
$\mu={\mathbb E}[-\log W_1],~\sigma^2={\rm Var}[-\log W_1]$.
If 
\eq\label{aslog}
\log\left[\sum_{j=1}^k \rho_j\right]=o(k)\,,~~{\rm as}~k\to\infty,
\en
then the strong law of large numbers holds, i.e. $K_n\sim \mu^{-1}\log n$ {\rm \,a.s.}. Moreover, the random variable
 $(K_n-{\mathbb E}\,[K_n])/\sqrt{{\rm Var\,}[K_n]}$ converges in law to the standard Gaussian distribution, whereas the
moments satisfy
\eq\label{mome}
{\mathbb E}[K_n]\sim {\log n\over \mu}\,,~~~{\rm Var}[K_n]\sim {\log n\over \sigma^2\mu^{-3}}\,.
\en
\end{proposition}
\begin{proof}
By the construction of $(U_j)\giv_{\rho}(H_k)$, we have a dichotomy: 
$U_n\in \,]H_{k},H_{k-1}]$ implies that either $U_n$ will enter the transformed sequence
or will get replaced by some $H_i\geq H_{k}$. 
Let $U_{1n}<\ldots <U_{nn}$ be the order statistics of $U_1,\ldots,U_n$.
It follows that 
\begin{itemize}
\item[(i)]   if $U_{jn}>H_k$ then $K_n\leq j+k\,$,
%\sum_{i=1}^k\rho_j$,
\item[(ii)] if $U_{mn}<H_k$ for $m=\sum_{i=1}^k \rho_k$ then  $K_n\geq k$.
\end{itemize}

\par  Define $\xi_n$ by
$U_{\xi_n,n}<1/n<U_{\xi_{n+1},n}$ and recall that $J_n$ was defined by 
 $H_{J_n+1}\leq 1/n<H_{J_n}$, thus
$\xi_n$ is the number of uniforms to the left of $1/n$, and $J_n$ follows the CLT.
Clearly, $J_n$ and $\xi_n$ are independent and $\xi_n$ is binomial$(n,1/n)$.
By (i), we have $K_n\leq J_n+\xi_n$
%\sum_{j=1}^{\lfloor\log n\rfloor}\rho_j$ 
where $\xi_n$ is 
approximately Poisson$(1)$, which yields 
%together with (\ref{aslog})
the desired upper bound.

\par  The lower bound is more delicate. Introduce 
$$\psi_n=c\sum_{j=1}^{\lfloor \log n\rfloor}\rho_j$$
where $c$ should be selected sufficiently large. Then  by the assumption 
(\ref{aslog}) $\log \psi_n=o(\log n)$, which is enough to assure that the number, say $L_n$, of the $H_k$'s
larger than $\psi_n/n$  is still asymptotic to $J_n$. 
Because $L_n$ is close to Gaussian with moments as in (\ref{mome}), an easy large deviation estimate
implies that the inequality 
$L_n<(c/2)\log n$ holds with probability at least $1-n^{-2}$.
On the other hand, the number of uniforms smaller $\phi_n/n$ is 
also close to Gaussian with both central moments about $\phi_n$, hence,
in view of $\phi_n>c\log n$, a similar estimate shows that this number is at least $(c/2)\log n$
with probability at least $1-n^{-2}$. 
By application of (ii) with $k=L_n$ we see that the lower bound $K_n>L_n$ holds 
up to an event of probability $O(n^{-2})$.
This completes the proof of the CLT. 
Finally, since  both $J_n$ and $L_n$ are 
asymptotic to $\mu^{-1}\log n$ almost surely, the Borel-Cantelli lemma
implies that the same is valid for $K_n$.
\end{proof}

\section{A continuous time process}

The sequential construction of $\Pi$ from the frequencies $(P_k)$ can be embedded in continuous time by 
letting the elements $1,2,\ldots$ 
arrive at epochs of a rate-$1$ Poisson process on ${\mathbb R}_+$. 
Let $R_t$ be the total frequency of the blocks which are not represented by 
the elements arrived before $t$, then $R=(R_t)$
is a nonincreasing 
pure-jump process with piecewise-constant paths and $R_0=1$.

\par Suppose as in Section \ref{s:CLT} that  $W_k$'s are independent and identically distributed. The process
$R$ is then easy to describe:
if after the $(k-1)$st  jump the process $R$ is in  state  $s$ then 
the time  in  this state
has distribution gamma$(\rho_k,s)$, and thereafter the state  is changed to $sW_k$.
The sojourns in consequtive states $1, W_1, W_1W_2,\ldots$ are independent.
The instance $\rho=(1,1,\ldots)$ corresponds to a known self-similar Markov process
which appears as a `tagged particle' process
in random fragmentation models \cite{BertoinBook}. 
For general $\rho$ the process is no longer
Markovian, as one needs to  also include the time spent in the current state  to 
summarise the history. 

\par A minor adjustment of Proposition \ref{CLT} to the continuous-time setting 
allows to conclude that under the same assumptions
the number of jumps during the time $[0,T]$ 
is approximately Gaussian, as $T\to\infty$.
In fact, the process $R$ is well defined for arbitrary  positive values $\rho_k$ ($k=1,2,\ldots)$,
in which case an analogous CLT is readily acquired  by interpolation from the case of integer $\rho_k$'s.

\section{Example: the chain records}

Next is an example of Pitman's partially exchangeable partitions, so the constraint is $\rho=(1,1,\ldots)$.
Consider a Borel space
${\cal Z}$ 
  endowed with a distribution $\mu$ and
some 
measurable strict partial order $\prec$.
For a sample $V_1,V_2,\ldots$
 from $({\cal Z},\mu)$, 
we say that a {\it chain record} occurs at index $j$ if either $j=1$, or $j>1$ and 
$V_j$ is $\prec$-smaller than the last chain record in the sequence $V_1,\ldots,V_{j-1}$.
The instance of ${\mathbb R}^d$ with the natural coordinate-wise partial order  was discussed in \cite{chain}.

\par Let $R_k,~k=1,2,\ldots$, be the sample values when the chain records occur; the sequence $(R_k)$ is a `greedy'
falling chain of the partially ordered sample $(V_j)$.
Introducing the lower sets
$L_v:=\{u\in {\cal Z}:u\prec v\}$, we define $D_k:=L_{R_k}$, $G_k:=D_{k}\setminus D_{k-1}$ (where $D_0:=\varnothing$), 
and we define a constrained 
exchangeable partition $\Pi=(B_k)$ as in Section \ref{intro}.
The frequencies of $B_k$'s are $P_k=\mu(L_{G_k})$, and we have $H_k=\mu(L_{R_k})$, as is easily seen.

\par To guarantee a `stick-breaking' form of $(H_k)$, as in Section \ref{6}, we need to  
assume a self-similarity property of the sampling space.
We may call $({\cal Z},\mu,\prec)$ {\it regenerative} if (i) $\mu(L_v)>0$
 for $\mu$-almost all points $v\in {\cal Z}$, and (ii) 
the lower section $L_v$ with conditional measure $\mu(\cdot)/\mu(L_v)$
is isomorphic, as a partially ordered probability space, to the whole space
$({\cal Z},\mu,\prec)$. Since all $L_v$'s are in this sense the same,
the $H_k$'s undergo stick-breaking with i.i.d. residual fractions whose distribution is the same as that of 
 $ L_{V_1}$.
Under the hypothesis of Proposition \ref{CLT}, 
the number of chain records among the first $n$ sample points is approximately Gaussian, since this number coincides with the number
of blocks of $\Pi_n$.
A class of regenerative spaces is comprised of the
 Bollob{\'a}s-Brightwell box-spaces \cite{BB},
which have  all
intervals
$\{u: v\prec u\prec w\}$ for $v\prec w$
 isomorphic to the whole space (and not only lower sections).
\par Further examples of regenerative spaces appear, in a disguise,
 in the context of multidimensional data structures like quad-trees or simplex-trees \cite{Devroye}.
More generally, constrained exchangeability  appears in connection with data structures 
which allow  multiple key storage at a node of the search tree.

\end{document}